\documentclass{article}
\usepackage{spconf,amsmath,graphicx}

\newcommand{\Sc}[1]{{\mathcal{#1}}}

\DeclareMathOperator*{\argmin}{argmin}

\title{Sparse Seismic Imaging using Variable Projection}
\name{Aleksandr Y. Aravkin$^{*\dagger}$, Tristan van Leeuwen$^*$, and Ning Tu$^*$\thanks{This work was in part financially supported by the
    Natural Sciences and Engineering Research Council of Canada
    Discovery Grant (22R81254) and the Collaborative Research and
    Development Grant DNOISE II (375142-08). This research was carried
    out as part of the SINBAD II project with support from the
    following organizations: BG Group, BPG, BP, Chevron, Conoco
    Phillips, Petrobras, PGS, Total SA, and WesternGeco.}}
\address{$^\dagger$ Computer Science, University of British Columbia\\ $^*$ Earth and Ocean Sciences, University of British columbia }
\begin{document}
%
\maketitle
\begin{abstract}
We consider an important class of signal processing problems where the signal of interest 
is known to be sparse, and can be recovered from data given auxiliary information about how 
this data was generated. For example, a sparse green's function may be recovered from 
seismic experimental data using sparsity optimization when the source signature is known. 
Unfortunately, in practice this information is often missing, and must be recovered from data
along with the signal using deconvolution techniques. 

In this paper, we present a novel methodology to simultaneously solve for the sparse signal and  
auxiliary parameters using a recently proposed {\it variable projection} technique. Our main contribution
is to combine variable projection with sparsity promoting optimization, obtaining an efficient algorithm 
for large-scale sparse deconvolution problems. We demonstrate the algorithm on a seismic imaging example. 

\end{abstract}
\section{Introduction}
\label{sec:intro}

Sparse regularization has proven to be an indispensable tool 
in many areas, including inverse problems~\cite{HerrmannJPHA2008} 
and compressive sensing~\cite{Donoho2006comp,candes06fdc}.  
If a signal $y$ is known to have a sparse or compressible (quickly decaying)
representation $y = Sx$, this information can be used to formulate optimization 
problems of the form 
\begin{equation}\label{BPDN}
\min \|x\|_1 \quad \text{s.t.} \quad \|ASx - b\|_2^2 \leq \sigma^2,
\end{equation}
where $A$ is a measurement matrix used to measure the true signal $y$, 
$b$ is a vector of data, and $\sigma$ is a threshold that depends on the 
characteristics of measurement error. In compressive sensing, 
it is possible to obtain recovery guarantees given properties of the true signal and $A$.
In more general inverse problems, these guarantees have not been found; 
it is therefore appropriate to consider~\eqref{BPDN} as a regularization approach to the 
least squares problem.
For example, in the seismic setting, where~\eqref{BPDN} has been
particularly useful~\cite{aravkin2012ICASSPfastseis}, $A$ is a linearized Born-scattering 
operator, $S$ is the curvelet transform, and $b$ is 
seismic data. While there are several popular algorithms that solve~\eqref{BPDN}, 
e.g. SPArsa~\cite{Wright2009}, the SPG$\ell_1$~\cite{BergFriedlander:2008} algorithm
has been particularly useful for seismic imaging~\cite{Herrmann2012a,herrmann11GPelsqIm,Li11TRfrfwi}.   

Many inverse problems contain unknown nuisance parameters that 
must be estimated in order to recover the solution~\cite{AravkinVanLeeuwen2012}. 
In seismic imaging, the source wavelet is typically unknown. 
The main contribution of this paper is 
to extend the approach of~\cite{AravkinVanLeeuwen2012} to the sparse inversion context, 
and derive simple modifications of standard sparse solvers to 
incorporate solutions of unknown nuisance parameters on the fly.

The paper proceeds as follows. In section~\ref{sec:OptForm}, we introduce
the seismic imaging problem with unknown wavelet, and formulate it as an 
extended sparsity promoting optimization problem~\eqref{EBPDN}. 
In section~\ref{sec:VarProj},
we review the ideas recently proposed in~\cite{AravkinVanLeeuwen2012}
that allow nuisance parameters to be estimated on the fly, and show how 
to incorporate these ideas into existing sparsity promoting formulations.   
In section~\ref{sec:ProjReg}, we develop an extended SPG$\ell_1$
algorithm to solve~\eqref{EBPDN}, 
and we present numerical results in section~\ref{sec:Numerics}.

\section{Imaging with Unknown Wavelet}
\label{sec:OptForm}

Seismic imaging is an approach to obtain a gridded subsurface velocity 
perturbation $y$ from seismic data, given a smooth starting model. 
Experiments are conducted by placing explosive
sources on the surface and recording the reflected waves with an array of
receivers on the surface. The data, $d_{i}$, in this case represents the Fourier
transform of the recorded time series for frequency $i$.  
The corresponding modeling operator, $F_{i}$, defines a linear relation between the recorded
data for the $i^{\mathrm{th}}$ frequency and the velocity perturbation.
The statistical model for data given $y$ is 
\begin{equation}\label{model}
d_i = \alpha_i F_i y + \epsilon_i,
\end{equation}
where $\epsilon_i$ is a statistical model for the measurement error, 
which is typically modeled as Gaussian, and $\alpha_i$ are unknown
complex source wavelet coefficients.  
Note that the model~\eqref{model} is no longer linear in the decision 
variables $(x,\alpha)$---it is bilinear.  
Since the perturbation $y$ is known to be sparse 
in the Curvelet frame $C$, formulation~\eqref{BPDN} has been 
successfully used to recover $y = Cx$~\cite{aravkin2012ICASSPfastseis}
when the source wavelet is known. In full generality, the 
joint inverse problem for the perturbation $y$ and wavelet $\alpha$
is given by 
\begin{equation}\label{EBPDN}
 \min_{x,\alpha} \|x\|_1 \quad \text{s.t} \quad 
 \sum_{i} \|d_i - \alpha_i F_i C x\|_2^2 \leq \sigma^2.
\end{equation}
Note that the $\alpha$ parameters make the problem
more difficult, because the forward model~\eqref{model}
is no longer linear in the decision variables $(x,\alpha)$,
and the problem~\eqref{EBPDN} is nonconvex.

\section{Variable Projection}
\label{sec:VarProj}

We begin 
by considering the problem
\begin{equation}\label{ELASSO}
\min_{x, \alpha}  \sum_{i} \|d_i - \alpha_i F_i C x\|_2^2 
\quad\text{s.t.}\quad \|x\|_1 \leq \tau\;.
\end{equation}
The relationship between~\eqref{ELASSO} and~\eqref{EBPDN}
will be fully explained in section~\ref{sec:ProjReg}.
In this section, we show how to use results from~\cite{AravkinVanLeeuwen2012}
to design an effective algorithm for~\eqref{ELASSO}.

If we define $\Sc{X} = \{x: \|x\|_1 \leq \tau\}$, 
problem~\eqref{ELASSO} is of the form
\begin{equation}
\label{InverseProblemClass}
\Sc{P} \quad \min_{x \in \Sc{X}, \alpha} g(x,\alpha)\;,
\end{equation}
where 
for any given $x \in \Sc{X}$, one can easily find 
\begin{equation}
\label{Projection}
\bar\alpha(x) = \argmin_{\alpha} g (x, \alpha)\;.
\end{equation}
In fact, $\bar\alpha(x)$ is available in closed form when 
the least squares penalty is used in~\eqref{ELASSO}.
The key idea in~\cite{AravkinVanLeeuwen2012}
is to consider the modified objective
\begin{equation}\label{projObj}
\tilde g(x) = g(x, \bar\alpha(x)), 
\end{equation}
using the convenient formula 
\begin{equation}\label{gradient}
\nabla_x \tilde g(\bar x) = \nabla_x g(\bar x, \bar\alpha(\bar x)).
\end{equation}
This is basically a generalization of the \emph{variable projection} algorithm
~\cite{golub03}.

In the current setting this means that instead of solving~\eqref{ELASSO},
we can simply solve the modified problem
\begin{equation}\label{ELASSOproj}
\min_{x}  \sum_{i} \|d_i - \bar\alpha_i(x) F_i C x\|_2^2 
\quad\text{s.t.}\quad \|x\|_1 \leq \tau
\end{equation}
using e.g. the projected gradient iteration
\[
x^{k+1} = P_{\Sc{X}}[x_k - \alpha_k \nabla_x \tilde g(\bar x^k)]
\]
with $ \nabla_x \tilde g$ computed via~\eqref{gradient},
with $\bar \alpha(x)$ given by~\eqref{Projection}.
By~\cite[Corollary 2.3]{AravkinVanLeeuwen2012}, a stationary point of~\eqref{ELASSOproj} 
is also stationary point of~\eqref{ELASSO}.

\section{Projected Regularized Inversion}
\label{sec:ProjReg}

In the previous section, we showed how to solve
the extended problem~\eqref{ELASSO}. However, 
the formulation~\eqref{EBPDN} is more important to us
from the modeling perspective, 
since it is always easier to provide a noise threshold
$\sigma$ than to figure out the `right' sparsity level $\tau$.
In fact, the SPG$\ell_1$ algorithm solves formulation~\eqref{BPDN} 
by solving a series of subproblems that find the sparsity level 
$\tau$ automatically given an input threshold $\sigma$. 

In this section, we extend this approach to the pair of problems~\eqref{EBPDN}
and~\eqref{ELASSO}. First, define
\begin{equation}\label{ValueFcn}
v(\tau) = \min_{x, \alpha}  \sum_{i} \|d_i - \alpha_i F_i C x\|_2^2 
\quad\text{s.t.}\quad \|x\|_1 \leq \tau
\end{equation}

Suppose we find $\bar \tau$ such that $v(\bar\tau) = \sigma^2$. 
Can we expect that the corresponding minimizers
of~\eqref{ELASSO} coincide with the minimizers
of~\eqref{EBPDN}?
This question is answered in surprising generality
by~\cite[Theorem 2.1]{AravkinBurkeFriedlander:2012}:
as long as any minimizer $\bar x$ 
of~\eqref{EBPDN} satisfies $\sum_{i} \|d_i - \alpha_i F_i C \bar x\|_2^2 = \sigma^2$, 
then the set of minimizers of $\eqref{EBPDN}$ and $\eqref{ELASSO}$ match, 
and $\|\bar x\|_1 = \bar \tau$ where $v(\bar \tau) = \sigma^2$. 

This general result points to using the following strategy: 
solve $v(\tau) = \sigma^2$ by Newton's method
\begin{equation}\label{NewtonIter}
\tau^{k+1} = \tau_k - \frac{v(\tau_k) - \sigma^2}{v'(\tau_k)}\;.
\end{equation}
This is in fact the strategy used by SPG$\ell_1$ to solve the 
problem~\eqref{BPDN}, for an appropriately defined value function.  
In order to implement this strategy, we have to be able
to evaluate both $v(\tau)$ and $v'(\tau)$ for~\eqref{ValueFcn}.

Evaluating $v(\tau)$ is straightforward: we simply use the projected
gradient method detailed in section~\ref{sec:VarProj}. 
However, $v'(\tau)$ is more difficult, since the most general 
variational results for value functions~\cite{AravkinBurkeFriedlander:2012}
require linearity of the forward model, which is violated by~\eqref{model}.
Nonetheless, if we treat $\bar \alpha_i := \bar \alpha_i(\bar x)$ as fixed, then by~\cite[Theorem 6.2]{AravkinBurkeFriedlander:2012}, we get 
\[
v'(\tau) \approx -\left\|\sum_i \bar \alpha_i C^TF_i^T(d_i - \bar\alpha_i F_iC\bar x)\right\|_{\infty}
\]
where $\bar x$ solves~\eqref{ELASSO} for $\tau$. 

\begin{figure}
\centering
\begin{tabular}{c}
\includegraphics[scale=.32]{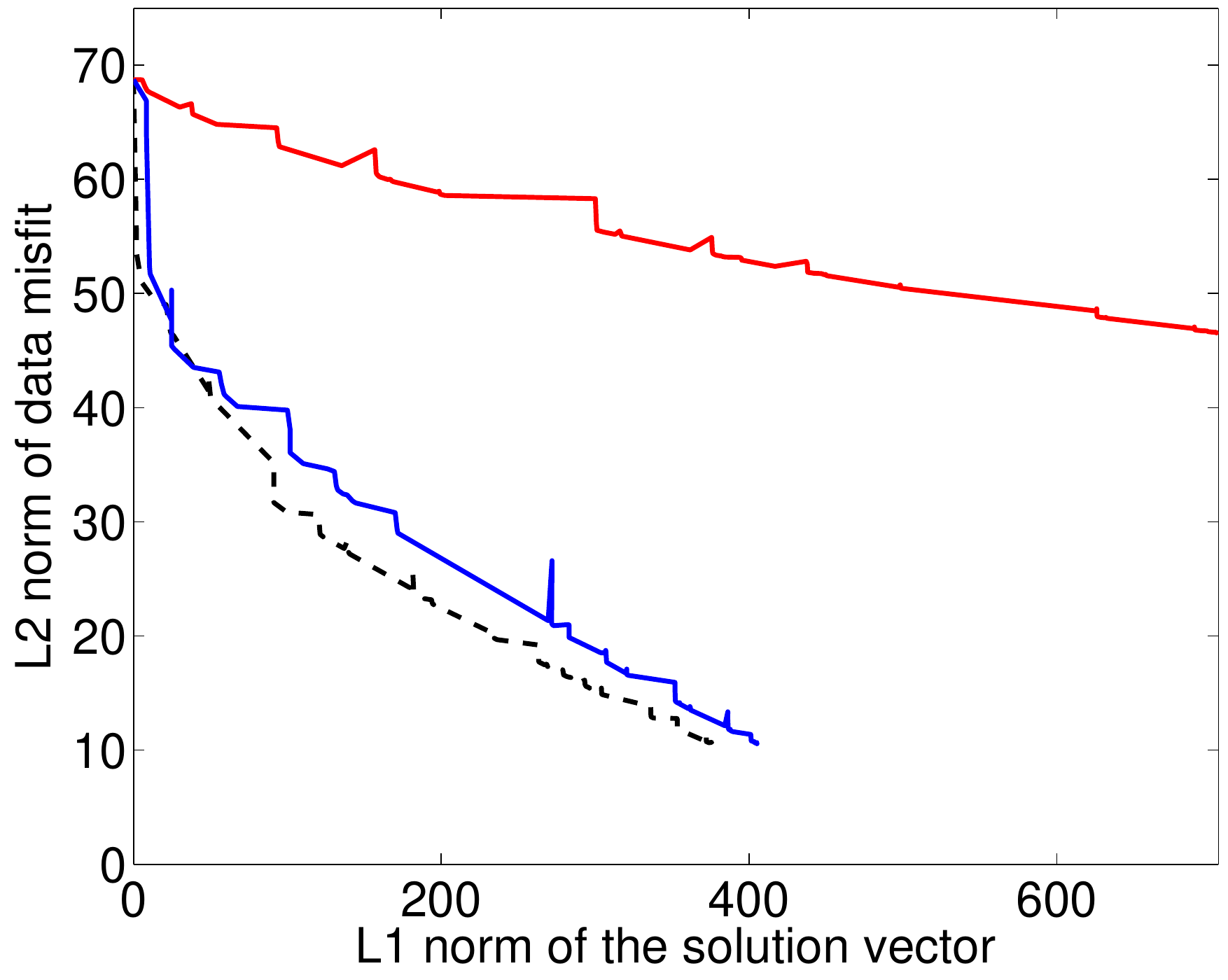}
{\small (a)}
\\
\includegraphics[scale=.32]{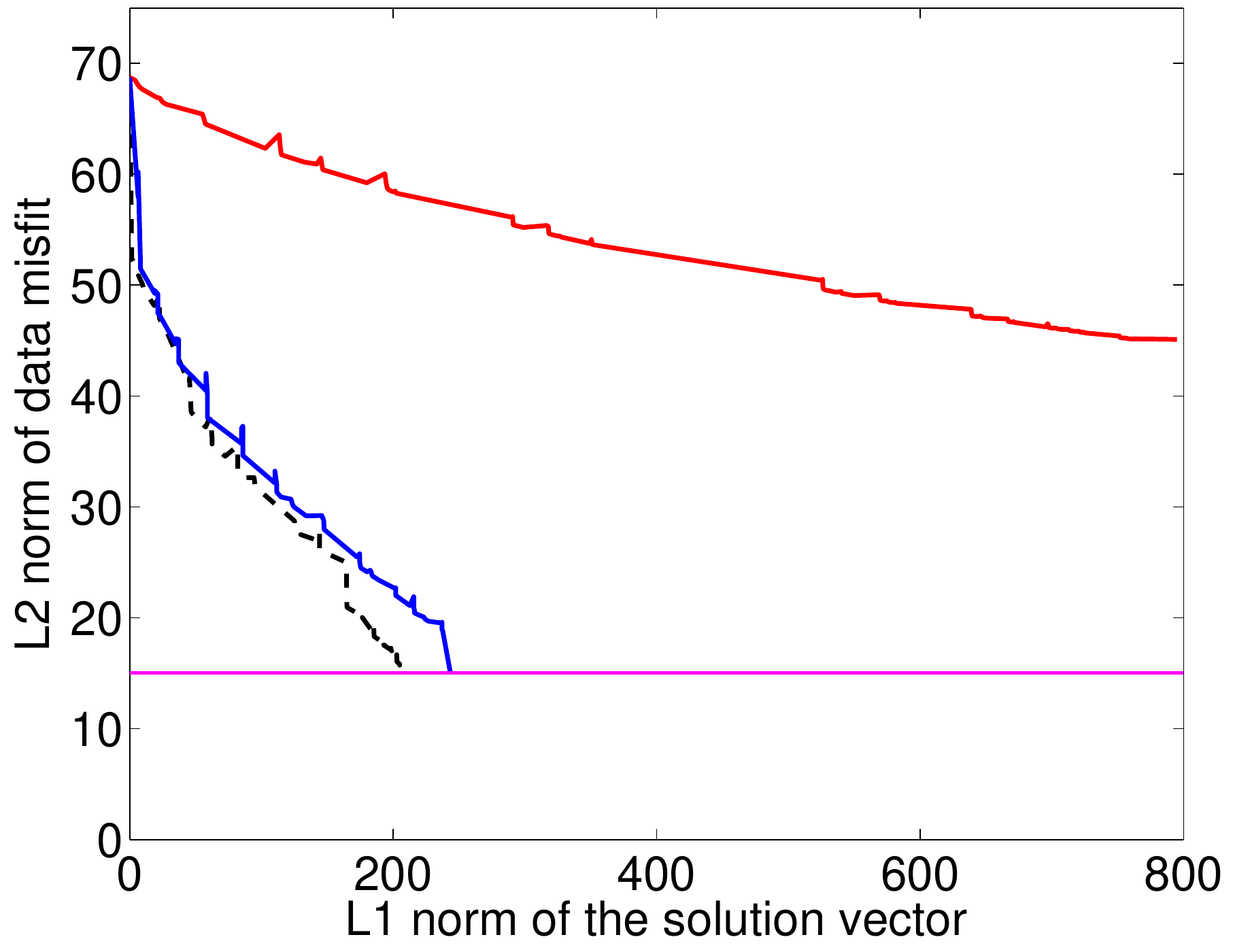}
{\small (b)}
\end{tabular}
\caption{Convergence histories using the true wavelet (dashed), a wrong source weight (red) and 
the estimated source weight (blue) when (a) $\sigma =0$ and (b) $\sigma = 15$.}
\label{fig:conv}
\end{figure}

We note that the expression above is an approximation to the derivative, 
and the quality of the approximation remains to be determined. 
If the source weight can be estimated fairly quickly (so that it is not changing 
significantly between iterations), the approximation above becomes exact. 
For the experiments in the next section, 
we found that the proposed Newton iteration gives 
nearly the same result as the one with a fixed, `true' source-weight. 
We also verified that when we pick $\sigma$ that is reachable 
within our computational budget of 150 iterations, 
the algorithm correctly finds the root $v(\tau) = \sigma$; 
see figure~\ref{fig:conv} (b). 

\section{Numerical Results}
\label{sec:Numerics}
For the experiments we use a Matlab framework for seismic imaging and modelling \cite{leeuwen12}, 
and the CurveLab toolbox \cite{curvelab}. Both of these are freely available for non-commercial purposes.
The algorithm to solve \ref{EBPDN} is based on the SPG$\ell_1$ code \cite{spgl1:2007}, which 
is also available for download. 

\begin{figure}
\centering
\begin{tabular}{c}
\includegraphics[scale=.35]{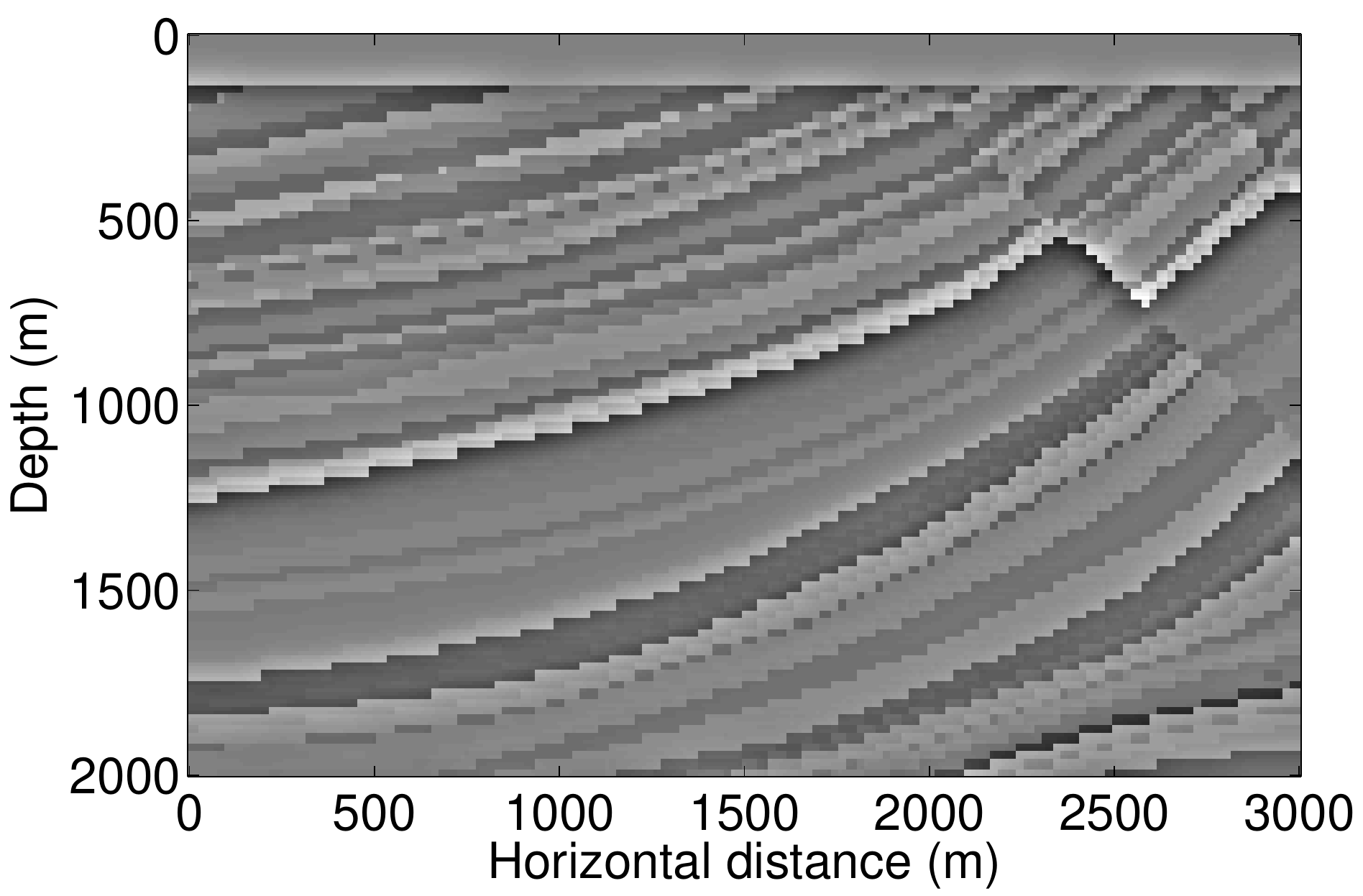}
\end{tabular}
\caption{True perturbation used for numerical experiment}
\label{fig:truemodel}
\end{figure}

We generate data for the velocity perturbation defined on a 201 x 301 grid with 10 m spacing
depicted in figure \ref{fig:truemodel} for 6 frequencies between 5 and 25 Hz, 
301 equispaced receivers and 15 composite sources, all located at the top of the model. 
We note that his leads to a underdetermined problem with 27090 equations and 60501 unknowns. 
We use SPG$\ell_1$ to solve (\ref{EBPDN}) either with $\alpha$ fixed or with $\alpha$ estimated
using the procedure outlined above. Since there is no noise in this example we use $\sigma = 0$
and run the algorithm for a fixed number of iterations (100 in this case).

Note that there is a fundamental non-uniqueness in the problem; if we multiply the source-weights
with a constant factor, we can compensate for this by dividing the reconstructed model
by the same factor. Therefore, we normalize the results such that the source-weights for each reconstruction
have the same norm (i.e., $\sum_i \bar\alpha_i(x)^2$ is the same for all reconstructions).

The reference result using the true source signature is shown in figure \ref{fig:reconmodel} (a).
If we do not estimate the source signature and use $\alpha_i = 1$, we do not get a good reconstruction, 
as is shown in figure \ref{fig:reconmodel} (b).
Finally, if we estimate the source signature according to the strategy outlined in this paper,
we obtain the result depicted in figure  \ref{fig:reconmodel} (c).

The (normalized) optimal source weights $\bar\alpha(x)$ evaluated at the final results as well as 
the true source weight are depicted in figure \ref{fig:wavelet}. 
This shows that our approach is able to recover both the model and the source weight.
The convergence histories of the SPG$\ell_1$ algorithm are shown in figure \ref{fig:conv}.

\section{Discussion and Conclusions}
\label{sec:Conclusions}

We have proposed a novel method for estimating nuisance 
parameters in the context of sparsity regularized inverse problems, 
and in particular we have focused on source wavelet estimation
in seismic imaging. 
The method draws on the idea of {\it variable projection}
in order to estimate nuisance parameters on the fly, 
and can be implemented via a straightforward
modification to existing sparse solvers. 

\begin{figure}
\centering
\begin{tabular}{c}
\includegraphics[scale=.35]{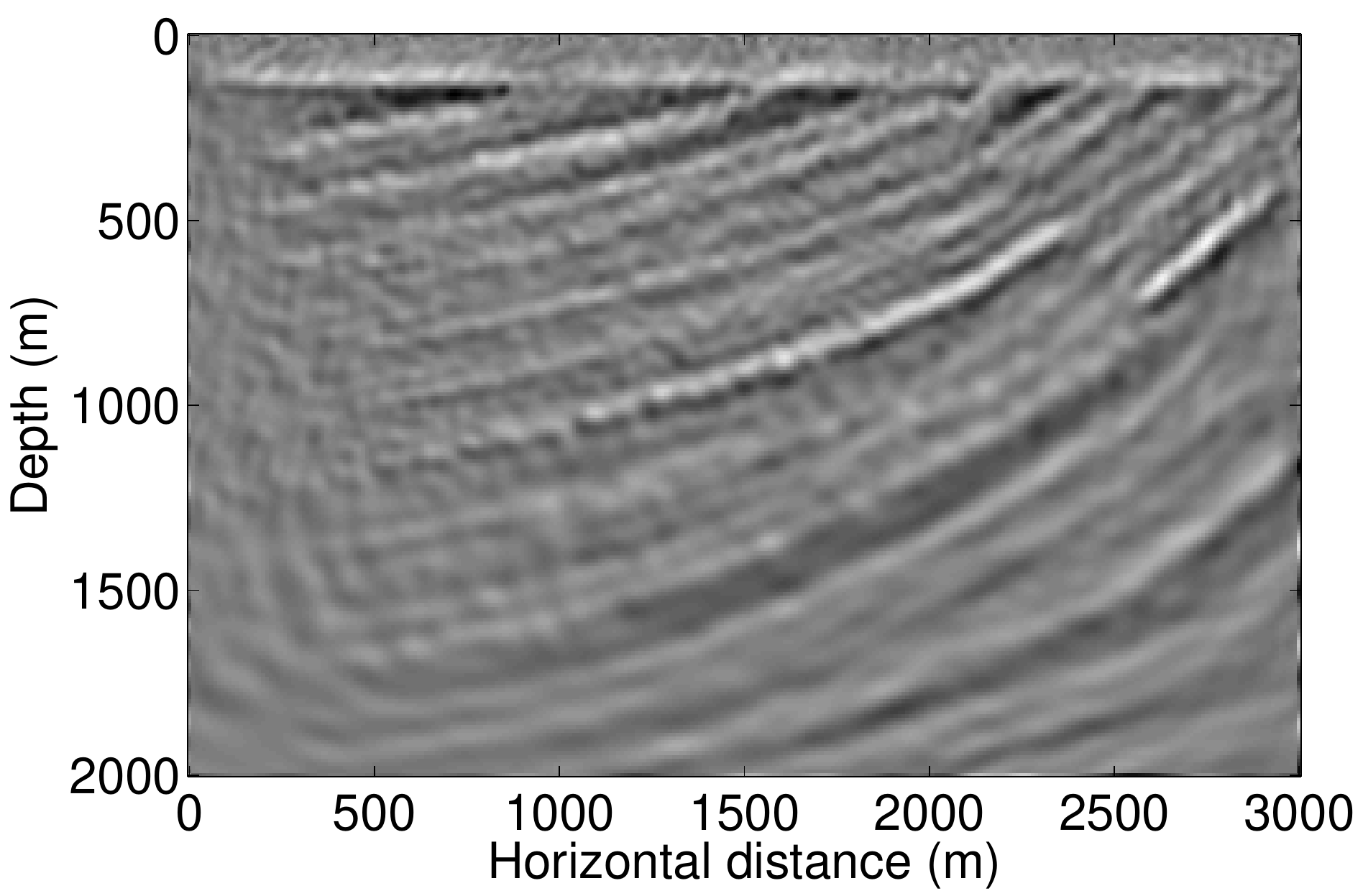}\\
{\small (a)}\\
\includegraphics[scale=.35]{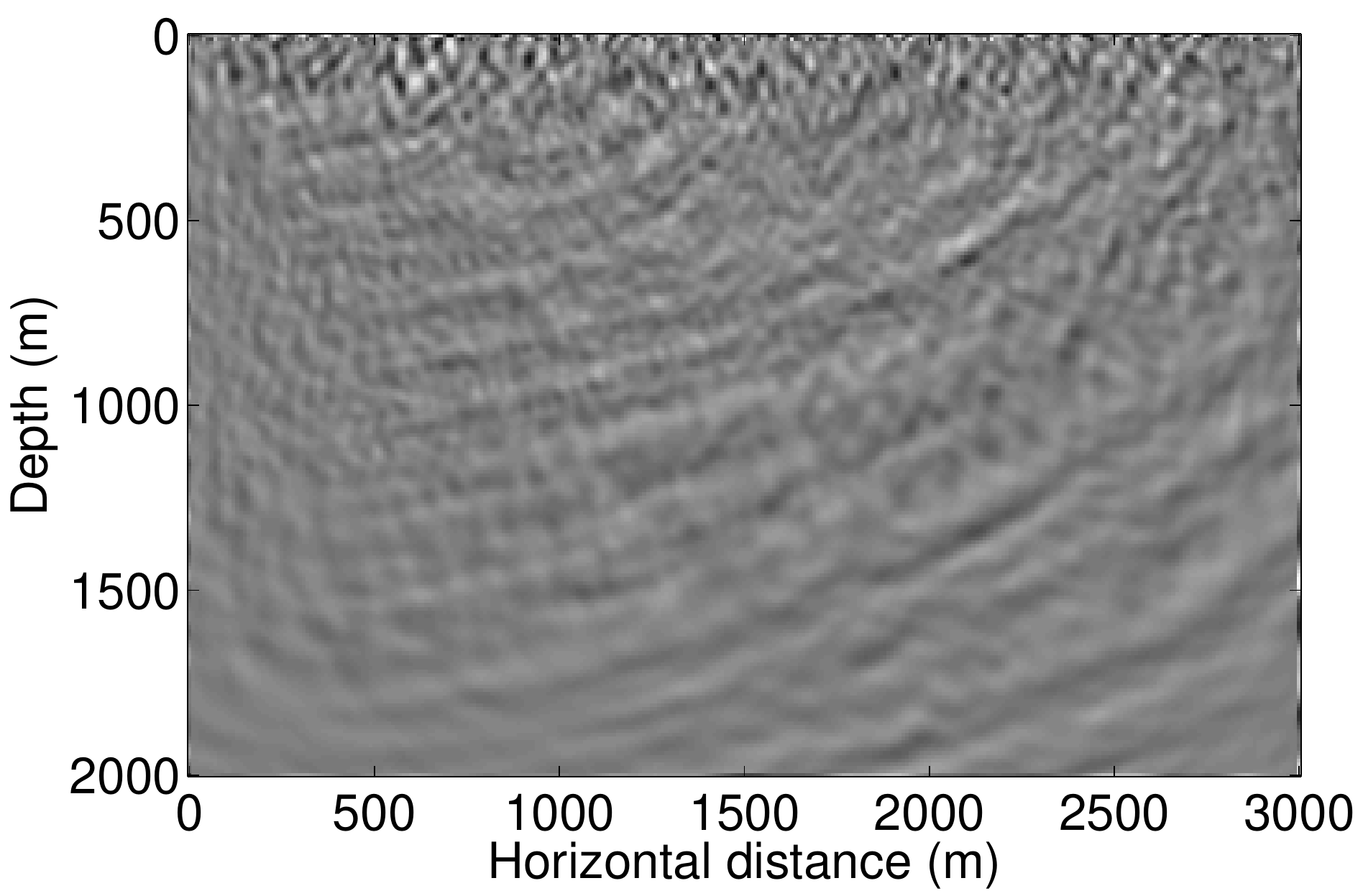}\\
{\small (b)}\\
\includegraphics[scale=.35]{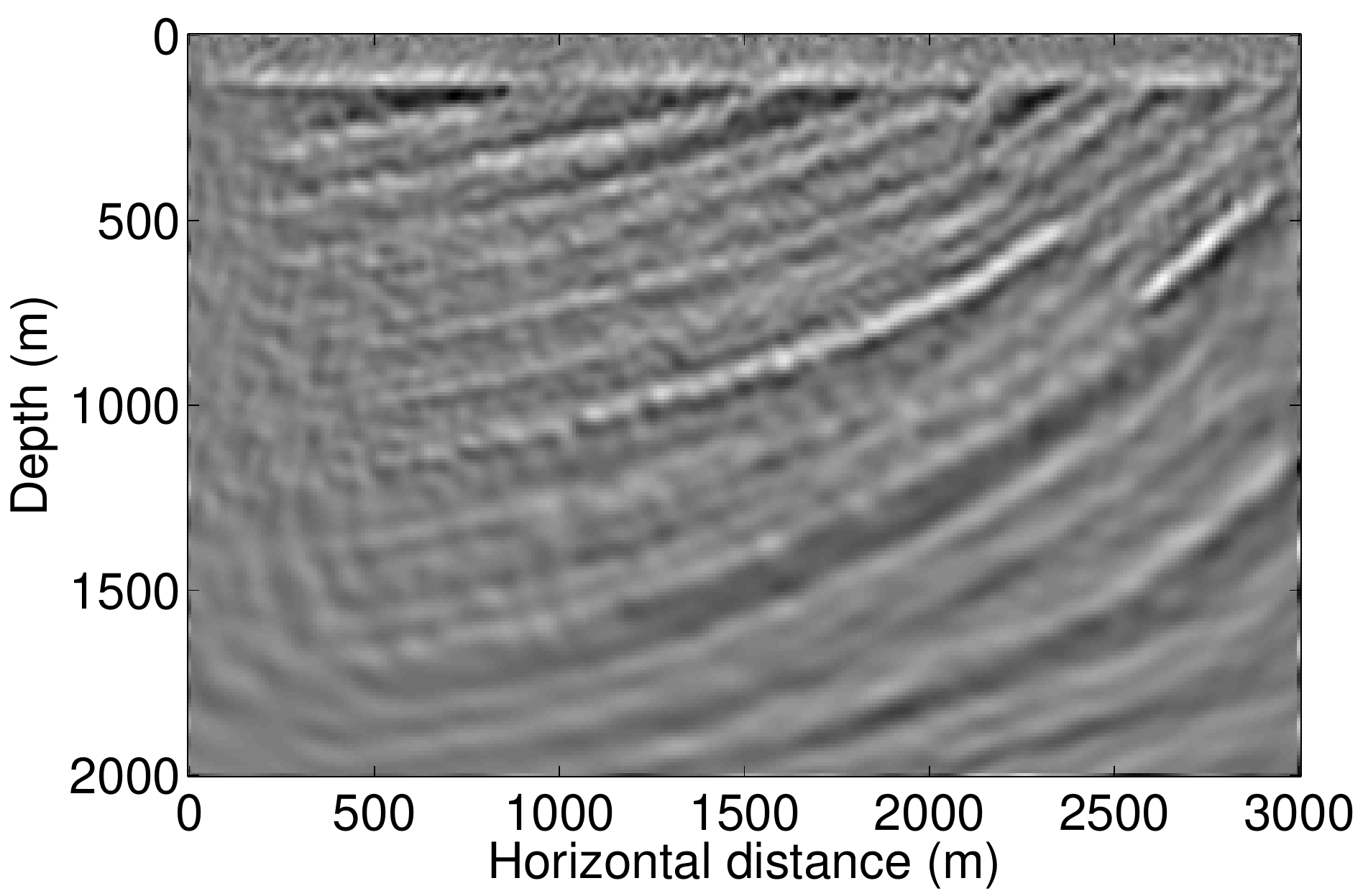}\\
{\small (c)}\\
\end{tabular}
\caption{Reconstructed models for the true wavelet (a), a wrong wavelet (b) 
and using the wavelet estimation procedure (c).}
\label{fig:reconmodel}
\end{figure}

\begin{figure}
\centering
\begin{tabular}{cc}
\includegraphics[scale=.32]{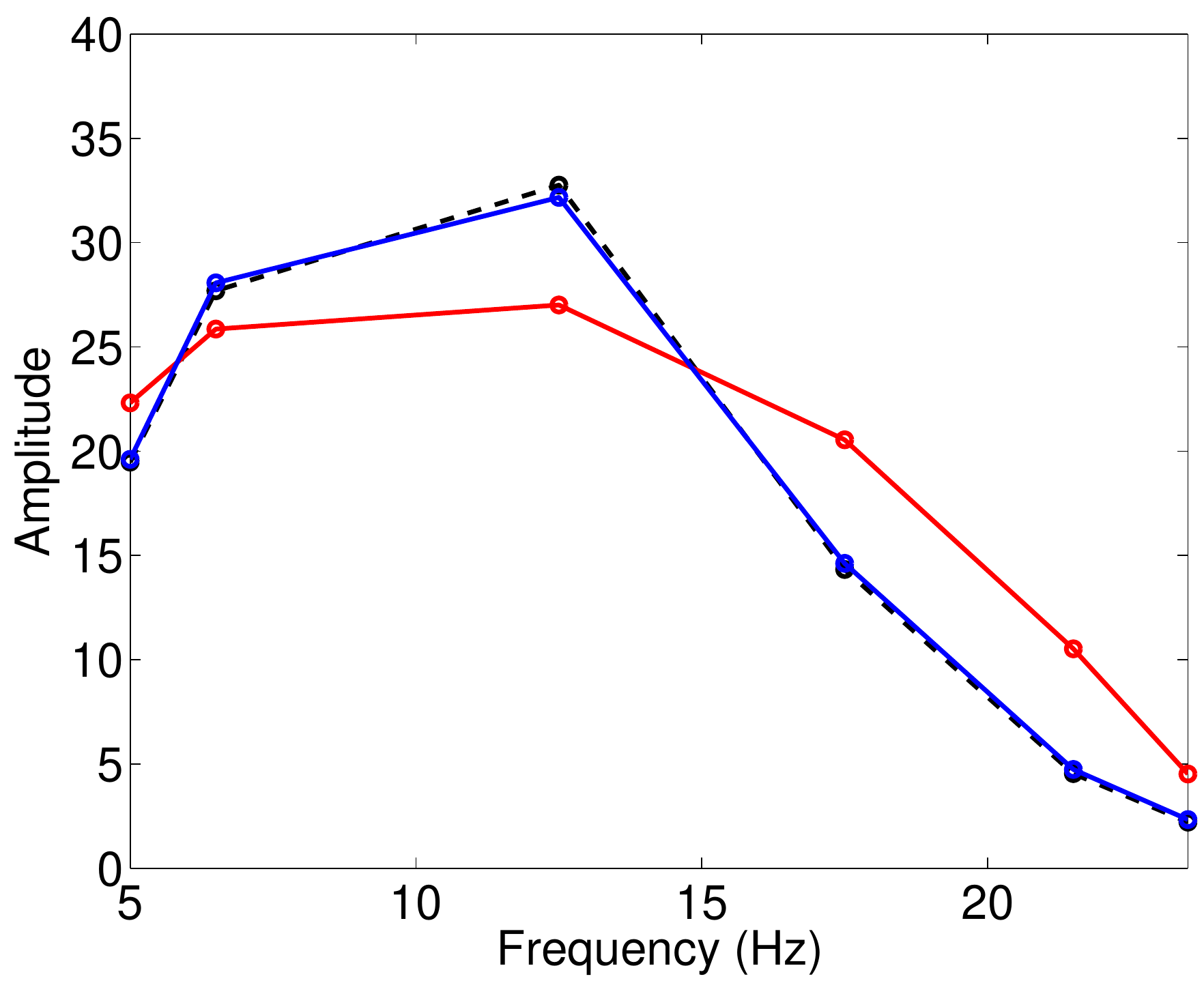}
\small(a)
\\
\includegraphics[scale=.32]{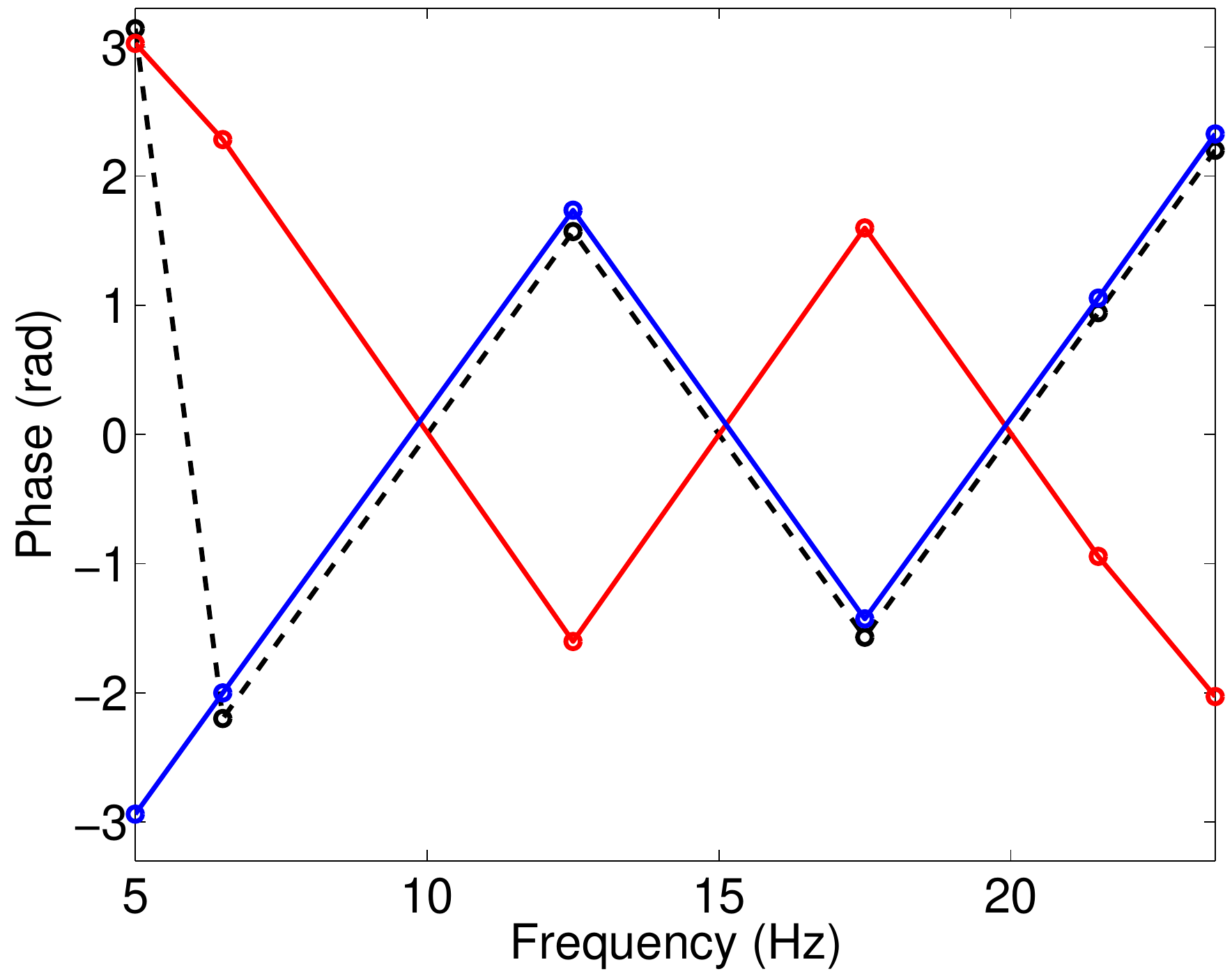}
{\small (b)}\\
\end{tabular}
\caption{Amplitude (a) and phase (b) of the optimal source weights evaluated at the 
models depicted in \ref{fig:reconmodel} (b) (red) and \ref{fig:reconmodel} (c) (blue). 
The true source weight is also shown (dashed line).}
\label{fig:wavelet}
\end{figure}

Numerical experiments demonstrate that 
the source wavelet can be recovered 
successfully in this manner (figure~\ref{fig:wavelet}),
and that the recovery of primary parameters
(specifically of the image) is improved when 
the wavelet is estimated (figure~\ref{fig:reconmodel}). 

Note that after section~\ref{sec:VarProj}, we can 
already solve~\eqref{ELASSO},  
but nonetheless lot of effort is devoted in section~\ref{sec:ProjReg}
to develop a method for solving~\eqref{EBPDN}.  
The main point here is that while it is difficult 
to come up with a reasonable value for $\tau$ in~\eqref{ELASSO},
it is straightforward to come up with a good 
value for $\sigma^2$ in~\eqref{EBPDN}. 
In fact, given a finite computational budget, and no
estimate for $\sigma^2$, one can always pick $\sigma^2 = 0$
and perform a fixed number of iterations. In this mode, the algorithm in
section~\ref{sec:ProjReg} solves several~\eqref{ELASSO}
problems inexactly, picking the corresponding sequence of $\tau$
values according to iteration~\eqref{NewtonIter}.
This is exactly what was done to obtain the numerical examples
in section~\ref{sec:Numerics}.

\clearpage
\bibliographystyle{IEEEbib}
\bibliography{mybib}

\end{document}